**Author for correspondence:**
Andreas Morr
e-mail: andreas.morr@tum.de

# Red noise in continuous-time stochastic modelling

Andreas Morr[1,2], Dörte Kreher[3] and Niklas Boers[1,2,4]

[1]Department of Aerospace and Geodesy, TUM School of Engineering and Design, Munich, Bavaria, Germany
[2]Department of Complexity Science, Potsdam Institute for Climate Impact Research (PIK) e V, Potsdam, Brandenburg, Germany
[3]Department of Mathematics, Humboldt University of Berlin, Berlin, Germany
[4]Department of Mathematics, University of Exeter, Exeter, UK

AM, 0000-0002-9804-5180; NB, 0000-0002-1239-9034

The concept of time-correlated noise is important to applied stochastic modelling. Nevertheless, there is no generally agreed-upon definition of the term red noise in continuous-time stochastic modelling settings. We present here a rigorous argumentation for the Ornstein–Uhlenbeck process integrated against time ($U_t\mathrm{d}t$) as a uniquely appropriate red noise implementation. We also identify the term $\mathrm{d}U_t$ as an erroneous formulation of red noise commonly found in the applied literature. To this end, we prove a theorem linking properties of the power spectral density (PSD) to classes of Itô-differentials. The commonly ascribed red noise attribute of a PSD decaying as $S(\omega) \sim \omega^{-2}$ restricts the range of possible Itô-differentials $\mathrm{d}Y_t = \alpha_t\mathrm{d}t + \beta_t\mathrm{d}W_t$. In particular, any such differential with continuous, square-integrable integrands must have a vanishing martingale part, i.e. $\mathrm{d}Y_t = \alpha_t\mathrm{d}t$ for almost all $t \geq 0$. We further point out that taking $(\alpha_t)_{t\geq 0}$ to be an Ornstein–Uhlenbeck process constitutes a uniquely relevant model choice due to its Gauss–Markov property. The erroneous use of the noise term $\mathrm{d}U_t$ as red noise and its consequences are discussed in two examples from the literature.

## 1. Introduction

In many fields of dynamical modelling, it is common practice to introduce a stochastic term to represent unresolved or uncertain dynamics within the otherwise deterministically defined differential equation [1–4]. To this effect, a general modelling approach for some one-dimensional observable $X$ is the stochastic differential equation

$$\mathrm{d}X_t = f(X_t, t)\mathrm{d}t + g(X_t, t)\mathrm{d}Y_t \tag{1.1}$$

for some suitable class of stochastic processes $Y$, to be specified later [5–7]. In the most classical setting, $Y$ will be the Wiener process $W$, leading to the so-called white noise model. This model exhibits independence in time, i.e. it is often said that the noise is $\delta$-correlated. However, in many cases, it is imperative to discard this assumption of independence, for instance, if the unresolved dynamics are suspected of exhibiting persistence in time [8,9]. Finding an appropriate model for such positively correlated noise will, in each application, depend on the observational or physical characteristics of the system [10,11]. One common model of correlated noise is so-called red noise [12–14]. The term is usually applied to mean a stochastic influence that has an exponentially decaying autocorrelation, and its frequency decomposition is dominated by low (red) frequencies [15–17]. There exists an established formulation of such a noise component in continuous-time stochastic modelling [18–23]. This is commonly denoted by $\mathrm{d}Y_t = U_t\mathrm{d}t$, where $U_t$ is an Ornstein–Uhlenbeck process. However, an alternative formulation has also been associated with the term red noise, namely $\mathrm{d}Y_t = \mathrm{d}U_t$ [24–27]. The two choices are not compatible, and only one satisfies common conceptions of red noise characteristics. We aim to clarify this confusion by introducing a fundamental result, constraining the possible noise models that can sensibly be termed red noise.

The field of stochastic analysis offers a wide range of possibilities for introducing stochasticity into continuous-time dynamics. For instance, the class of Itô-processes

$$Y_t = Y_0 + \int_0^t \alpha_s \mathrm{d}s + \int_0^t \beta_s \mathrm{d}W_s,$$

for suitable processes $\alpha = (\alpha_t)_{t\in\mathbb{R}_+}$ and $\beta = (\beta_t)_{t\in\mathbb{R}_+}$, allows for the interpretation of (1.1) as an equation of Itô-integrals. This work aims to constrain the class of such processes $Y$ for which it can be said that $\mathrm{d}Y$ constitutes red noise.

To achieve this, we will invoke properties of the power spectral density (PSD) that are commonly ascribed to the notion of red noise and introduce a new result about their implication on associated Itô-differentials. After finding that only differentials of the form $\mathrm{d}Y_t = \alpha_t\mathrm{d}t$ are admissible candidates for red noise, we show the unique relevance of choosing $\alpha$ to be an Ornstein–Uhlenbeck process as per its Gauss–Markov property. We clarify the common misrepresentation of red noise by $\mathrm{d}U$ along two examples from the literature in §3 and discuss the model restrictions we have made in §4.

## 2. Characterization through the power spectral density

Throughout this section, we will work on a filtered probability space $(\Omega, \mathcal{F}, (\mathcal{F}_t)_{t\in\mathbb{R}_+}, \mathbb{P})$ supporting a Brownian motion $W = (W_t)_{t\in\mathbb{R}_+}$. For any $Y_0 \in \mathcal{F}_0$ and predictable processes $\alpha = (\alpha_t)_{t\in\mathbb{R}_+}$ and $\beta = (\beta_t)_{t\in\mathbb{R}_+}$ the Itô-process

$$Y_t = Y_0 + \int_0^t \alpha_s \mathrm{d}s + \int_0^t \beta_s \mathrm{d}W_s, \ \ t \geq 0,$$

is well-defined if, e.g. the integrability condition

$$\int_0^T \mathbb{E}\left[\alpha_t^2 + \beta_t^2\right] \mathrm{d}t < \infty \quad \text{for all} \quad T > 0 \tag{2.1}$$

is satisfied. The PSD of $\mathrm{d}Y$ is defined as

$$S_{\mathrm{d}Y}(\omega) := \lim_{T\to\infty} \mathbb{E}\left[\frac{1}{T}\left|\int_0^T \exp(-i\omega t)\mathrm{d}Y_t\right|^2\right], \tag{2.2}$$

if the limit exists. One immediate result through the Itô-isometry is $S_{\mathrm{d}W} \equiv 1$ (see figure 1a). In this case, the power is evenly distributed across all frequencies in the spectral density, making the terminology 'white noise' a natural choice for this model.

If $\mathrm{d}Y_t = \alpha_t\mathrm{d}t$ for some stationary, centred and square-integrable process $\alpha$ with absolutely integrable autocovariance structure $R_\alpha(\tau)$, then the Wiener–Khinchin theorem applies:

$$S_{\alpha_t\mathrm{d}t}(\omega) = \mathcal{F}[R_\alpha(\tau)](\omega) := \int_{-\infty}^{\infty} \exp(-i\omega\tau)R_\alpha(\tau)\mathrm{d}\tau. \tag{2.3}$$

In most applications of correlated noise, the properties of the PSD are prominently featured. In fact, the name red noise stems from the observation that low frequencies exhibit the largest amplitudes in the PSD.

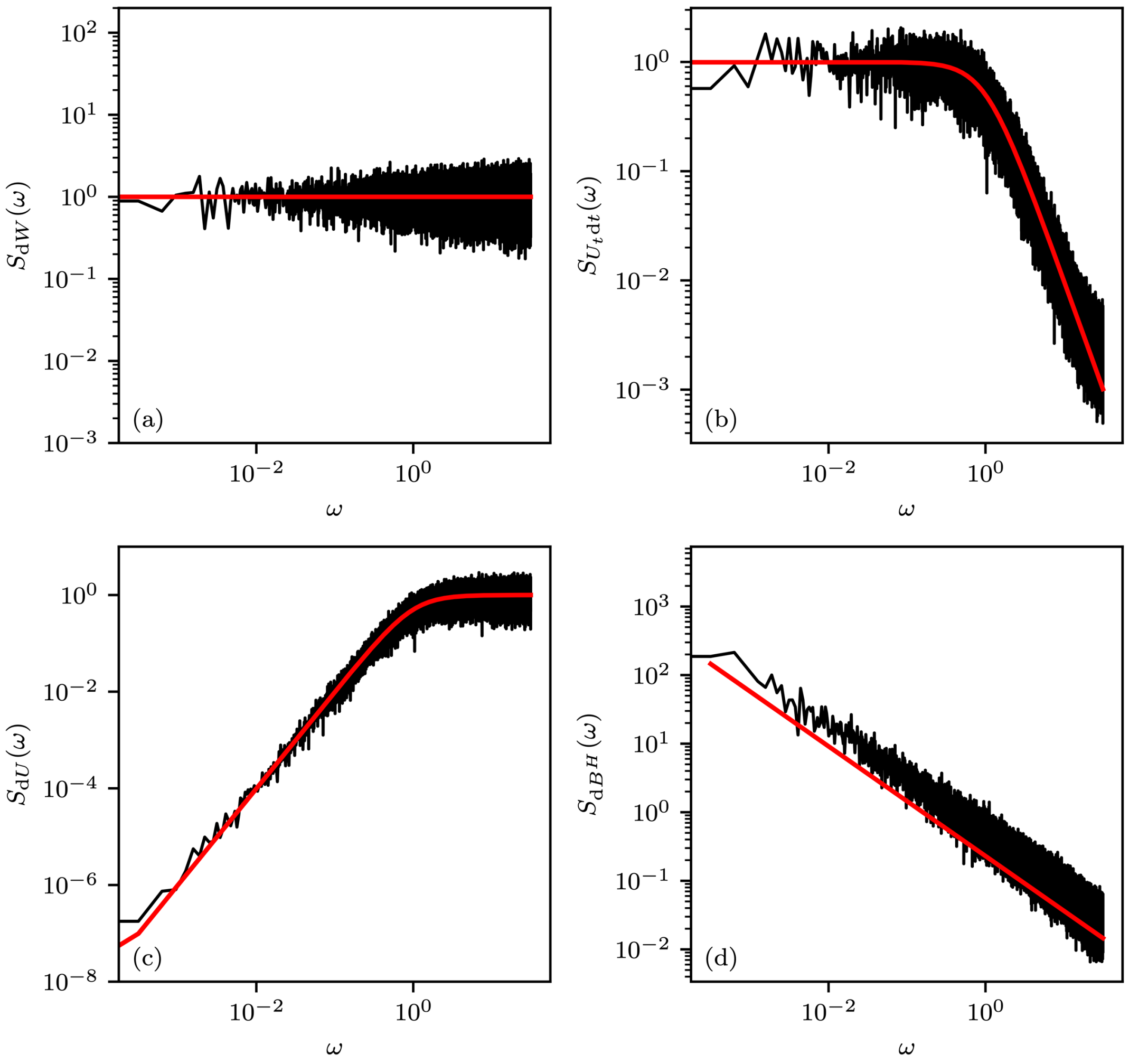

**Figure 1.** The theoretically computed (red) and observed (black) PSD for each of the noise differentials discussed in this work, shown on a log-log-scale. The observed PSD was obtained by taking the squared absolute value of the Fourier-transformed noise signal and subsequently averaging over 10 neighbouring frequencies, respectively. The length of the time-series samples is $T = 2 \cdot 10^4$ time units sampled at time steps $\Delta t = 10^{-1}$. (a) PSD of unit white noise $\mathrm{d}W$. (b) PSD of the red noise differential $U_t\mathrm{d}t$. The characteristic $\omega^{-2}$ asymptotic can be empirically observed until close to the Nyquist frequency, at which point the discrete nature of the data skews the PSD. (c) PSD of the differential $\mathrm{d}U$ sometimes erroneously referred to as red noise (see §3). (d) PSD of fractional Gaussian noise $\mathrm{d}B^H$, a red noise alternative outside the Itô-framework of integration (see §4). Parameter values were chosen as $\Theta = 1$ and $H = 0.9$.

A sufficient but not equivalent condition is to demand an asymptotically vanishing PSD, i.e. $S(\omega) \to 0$ as $\omega \to 0$. In much of the applied literature on this topic, a rate of decay of $\mathcal{O}(\omega^{-2})$ is taken to be the defining characteristic of red noise [13,14,28] and the noise instances are sometimes constructed directly through its PSD [15,16,29]. The specific dependence on $\omega$ is usually derived from observations in frequency space [30].

With this defining property of red noise, we are now able to formulate our main result, strongly constraining the possible choices of such continuous-time models in the form of an Itô-differential.

**Theorem 2.1.** *Let $\alpha$ and $\beta$ be adapted processes satisfying the integrability condition* (2.1). *Suppose that $\alpha$ has continuous paths and $\beta$ is predictable. Define the Itô-process*

$$Y_t = \int_0^t \alpha_s \mathrm{d}s + \int_0^t \beta_s \mathrm{d}W_s, \quad t \geq 0.$$

*1. Finite time horizon. Assume that the finite-time PSD of* $\mathrm{d}Y$ *vanishes in the limit of infinitely high frequencies, i.e.*

$$S_{\mathrm{d}Y}^{(T)}(\omega) := \mathbb{E}\left[\frac{1}{T}\left|\int_0^T \exp(-i\omega t)\mathrm{d}Y_t\right|^2\right] \xrightarrow[\omega\to\infty]{} 0.$$

*Then* $\beta_t = 0\,\mathbb{P}$*-a.s. for almost all* $t \in [0, T]$.

*2. Infinite time horizon. Assume that* $\alpha$ *and* $\beta$ *are stationary processes and let* $\alpha$ *be centred around* 0 *with an absolutely integrable autocovariance structure* $R_\alpha(\tau)$*. Assume that the PSD of* $\mathrm{d}Y$ *on an infinite time horizon exists and that it vanishes in the limit of infinitely high frequencies, i.e.*

$$S_{\mathrm{d}Y}(\omega) := \lim_{T\to\infty} \mathbb{E}\left[\frac{1}{T}\left|\int_0^T \exp(-i\omega t)\mathrm{d}Y_t\right|^2\right] \xrightarrow[\omega\to\infty]{} 0.$$

*Then* $\beta_t = 0\,\mathbb{P}$*-a.s. for all* $t \geq 0$.

We refer to appendix A.1 for a proof of this theorem and to appendix A.2 for an extension to the case with an additional fractional Gaussian noise term. If the objective of finding a suitable red noise Itô-differential $\mathrm{d}Y$ is to have it exhibit a vanishing PSD in the limit of infinitely high frequencies, then all choices necessarily have the form $\mathrm{d}Y_t = \alpha_t\mathrm{d}t$. The conditions on $\alpha$ and $\beta$ in the theorem are narrower if one examines an infinite time horizon since asymptotic behaviour needs to be taken into account. The strict stationarity condition for $\alpha$ and $\beta$ may be replaced by other suitable constraints on their asymptotic behaviour. In the finite time horizon case, the restrictions on $\alpha$ and $\beta$ reduce to being square-integrable together with the path-continuity of $\alpha$. If $\alpha$ and $\beta$ are themselves models of physical processes, these properties are often desired. For instance, Newtonian motion is continuous and ergodic measures of bounded chaotic motions are square-integrable.

Still, the form $\mathrm{d}Y_t = \alpha_t\mathrm{d}t$ allows for many options in the process $\alpha$. One specific choice comes to mind, as it is a common time-series model used for observed red noise [17,28]. The Ornstein–Uhlenbeck process $U_t$ following the Itô-diffusion

$$\mathrm{d}U_t = -\theta U_t\mathrm{d}t + \mathrm{d}W_t,\ U_0 \sim \mathcal{N}\left(0, \frac{1}{2\theta}\right) \tag{2.4}$$

has the autocovariance function (ACF) (see figure 2b).

$$R_U(\tau) = \frac{1}{2\theta}\exp(-\theta\tau). \tag{2.5}$$

Hence, $\theta$ determines the correlation decay of the process $U$. Equations (2.3) and (2.5) imply a PSD of (see figure 1b)

$$S_{U_t\mathrm{d}t}(\omega) = \mathcal{F}[R_U(\tau)](\omega) = \frac{1}{\theta^2 + \omega^2} = \mathcal{O}(\omega^{-2})$$

and in much of the applied literature, the term $\mathrm{d}Y_t = U_t\mathrm{d}t$ is indeed referred to as red noise [18–23].

Apart from the apt PSD of this model, we provide further conceptual motivation for this specific construction. First, since noise models are often physically motivated by homogenization and central limit theorem arguments, $\alpha$ will generally be a Gaussian process. If the homogenized dynamics are time-locally constant in nature, it is reasonable to look for a stationary $\alpha$. In most applications, it is desirable for the modelled red noise to exhibit the Markov property [28,31,32]. This is naturally the case if $\alpha$ is itself assumed to be the product of an Itô diffusion, i.e. $\mathrm{d}\alpha_t = \rho(\alpha_t)\mathrm{d}t + \sigma(\alpha_t)\mathrm{d}W_t$ with Lipschitz continuous $\rho$, $\sigma : \mathbb{R} \to \mathbb{R}$. We can leverage Theorem 1.1 in [33], which states that all stationary, measurable processes $\alpha$, which are simultaneously Gaussian and Markov, are of the Ornstein–Uhlenbeck type. Although the presupposition of the Markov property is not the most dominant characteristic of red noise posed in the literature, it is worth pointing out its consequences in this context.

Looking at this model through the lens of numerical implementation and discretization, the connection to the autoregressive process of order 1 (AR(1)) offers further motivation. Through the application of the Euler scheme, the noise term $dY_t = U_t dt$ becomes $U_{k\Delta t}\Delta t$. Sampling the process $U$ at evenly spaced time steps $\Delta t > 0$, one obtains the AR(1)-process

$$U_{(k+1)\Delta t} = \exp(-\theta\Delta t)U_{k\Delta t} + \frac{1}{\sqrt{2\theta}}\left(1 - \exp(-2\theta\Delta t)\right)^{1/2} z_k,$$

where the $(z_k)_{k\in\mathbb{N}}$ are i.i.d. unit Gaussian. The AR(1) process is the basis for most discrete-time stochastic models invoking the name red noise [13,25,34–37] and it too is a stationary Gauss–Markov process. If $U$ is itself an observable of some natural system, the correlation parameter $\theta$ can be statistically inferred on discrete-time [38,39] and continuous-time data [40,41]. When $U$ is introduced as temporal disturbances into some observed dynamics, there exist further avenues of indirect parameter estimation [42,43].

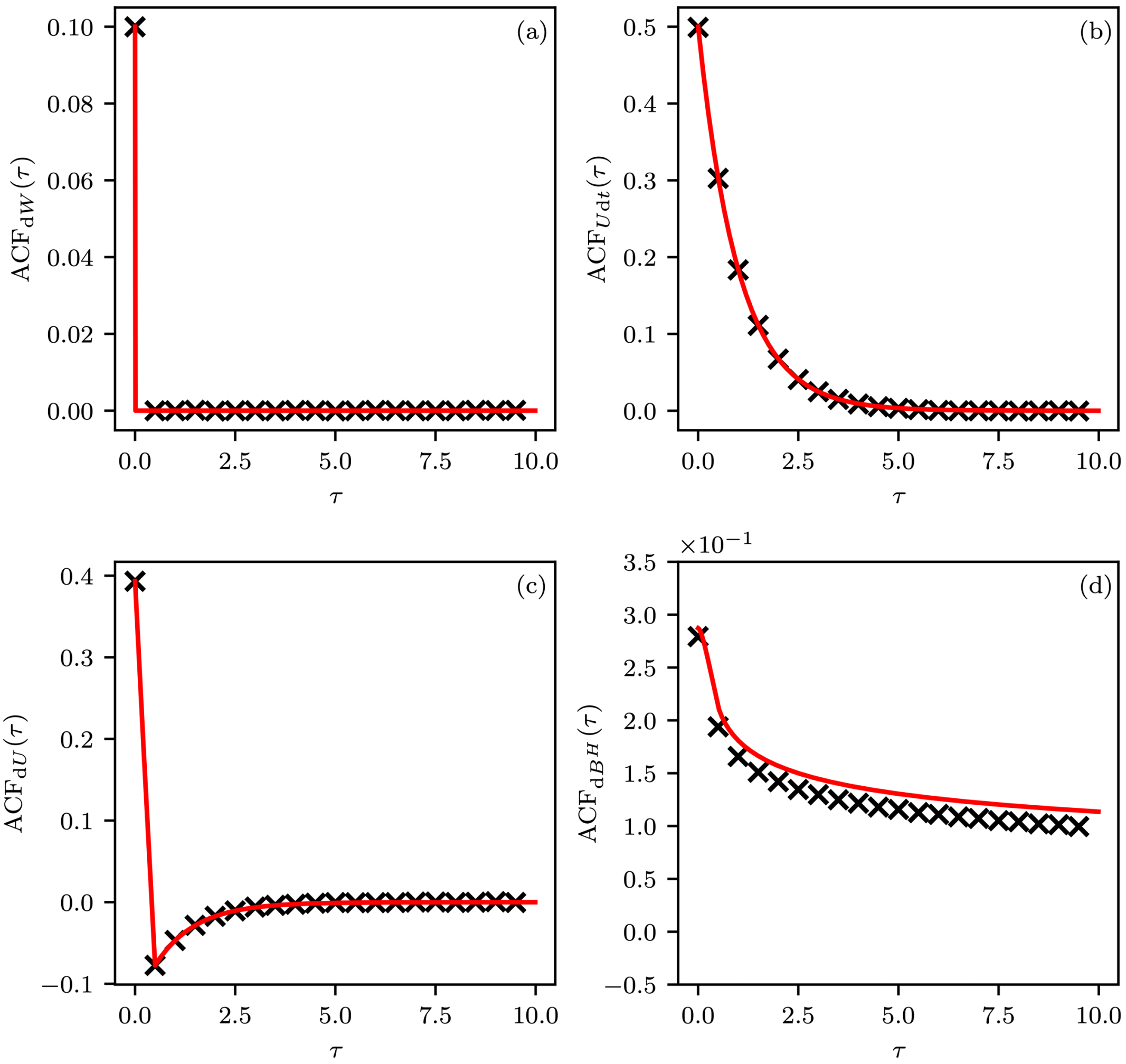

**Figure 2.** The theoretically computed (red) and observed (black) ACF for each of the discussed noise differentials. The length of the time-series samples is $T = 2 \cdot 10^5$ time units sampled at time steps $\Delta t = 0.5$. (a) ACF of unit white noise $dW$. (b) ACF of the red noise differential $U_t dt$, exhibiting a characteristic exponential decay of positive correlations. (c) ACF of the differential $dU$, a noise model with negative correlation in time. (d) ACF of fractional Gaussian noise $dB^H$, a red noise alternative outside the Itô-framework of integration (see §4). The latter noise model exhibits long memory with a correlation decay of order $\mathcal{O}(\tau^{2H-2})$ and the slight deviation of the observed ACF from theory is explained by the difficulty to accurately and efficiently generate samples of this process. Parameter values were chosen identically to the setting of figure 1.

We may conclude that the Ornstein–Uhlenbeck process $U$ defined by equation (2.4) constitutes a unique way of modelling red noise in continuous time through

$$dY_t = \alpha_t dt = U_t dt.$$

This noise model allows for an interpretation of equation (1.1) as a random ordinary differential equation:

$$\frac{dX_t}{dt} = f(X_t, t) + g(X_t, t)U_t.$$

There exists a natural limit of this red noise model to the canonical white noise model of uncorrelated disturbances. This occurs when $\theta \to \infty$ and considering $\mathrm{d}Y_t = \theta U_t \mathrm{d}t$ [18,44]. For this reason, another popular formulation of the Ornstein–Uhlenbeck SDE is

$$\mathrm{d}U_t = -\frac{1}{D} U_t \mathrm{d}t + \frac{1}{D} \mathrm{d}W_t,$$

for which a vanishing correlation time $D$ naturally yields the white noise model. The here discussed red noise model is commonly featured in the applied literature. However, no comprehensive justification using the breadth of its statistical characteristics and associated physical implications has hitherto been given in its favour.

# 3. Formulation of continuous-time noise models

In the applied literature, there exists another representation of continuous-time red noise, which we would like to discern here [24–27]. Instead of $\mathrm{d}Y_t = U_t \mathrm{d}t$, it is asserted that $\mathrm{d}Y_t = \mathrm{d}U_t = -\theta U_t \mathrm{d}t + \mathrm{d}W_t$, effectively choosing the time-derivative of $U$ as a noise model. It may seem visually appealing to label this red noise since $U$ constitutes a red noise process, and SDEs conventionally involve a differential. However, as implied by our Theorem 2.1, the PSD of $\mathrm{d}U$ exhibits the opposite of what we have so far understood as red noise characteristics (see appendix A.3 for a derivation):

$$S_{\mathrm{d}U} = \frac{\omega^2}{\theta^2 + \omega^2}.$$

For low frequencies, $S_{\mathrm{d}U}$ tends to 0, and it is monotonically increasing in $\omega$ (see figure 1c). The discrete-time noise terms resulting from discretizing such a differential through the Euler–Mayurama method at integration time-step $\Delta t$ would also be negatively correlated (here $\tau \geq \Delta t$) (see figure 2c):

$$\mathrm{Cov}\,(U_{\Delta t} - U_0, U_{\tau+\Delta t} - U_\tau) = \frac{1}{\theta} \exp(-\theta\tau)(1 - \cosh(\theta\Delta t)) < 0.$$

The usage of $\mathrm{d}U$ in this context constitutes a common misconception about the formulation of continuous-time stochastic models from desired discrete-time characteristics. If one posits a certain distribution or correlation in the noise component of a system and encounters the desired property in a stochastic process $V$, then the differential $\mathrm{d}V$ is, in general, not a suitable noise term since it may exhibit entirely different properties. Both the regularity and smoothness of sample paths, as well as the resulting statistical properties, will generally differ from the intended effect of the noise. In [24] and [25], this would instead imply using $V_t \mathrm{d}t$. Both of these studies posit an equivalence of a discrete-time system driven by AR(1) red noise and a continuous-time system driven by $\mathrm{d}U_t$ noise. This is an erroneous assertion which can result in faulty statistical models for physical systems, as we discuss further below. Even when introducing $\mathrm{d}U$ more generally as coloured noise [45,46], one should be aware of the conceptual implications of negative correlation. The two cited studies introduce $\mathrm{d}U$ as a stochastic term representing unresolved, fast dynamic variables in the context of electronics and climate (discussed in the following paragraph), respectively. If the aim of introducing correlation is to model temporal persistence in the noise forcing, then this would heuristically always call for positive correlation and $\mathrm{d}U$ would not be a good choice.

In each individual modelling task, there may still be reasons to choose either of the above noise models $U_t \mathrm{d}t$ or $\mathrm{d}U$. However, for many generic applications, the use of $\mathrm{d}U$ to model (negatively) correlated noise is not feasible given its non-dispersive property. To see this, we consider one setting in which the effect of correlated noise has recently gained attention [25,43,45,47–49]. In systems undergoing abrupt transitions by means of a changing stability landscape (bifurcation), the type of noise inherent to the system can have a substantial influence on whether the transition can be anticipated. This is because classical warning signals rooted in the phenomenon of critical slowing down (CSD) [43,50,51] are based on trends in time-series statistics, which are themselves dependent on the nature of the random disturbances. If the purpose of these disturbances is to model the system state's random dispersal, then the deterministic system's non-hyperbolic equilibria should not persist (see figure 3a–d and its caption). That is to say, once a deterministic equilibrium loses stability by means of a bifurcation and becomes non-hyperbolic, the system state should meander away in a manner akin to Brownian particle motion. Instead, systems forced by low amplitude noise of the form $\mathrm{d}U$ will maintain a quasi-stable state at the non-hyperbolic equilibrium and assume the statistics of $U$. Hence, the variance and autocorrelation of the system will remain bounded by the values associated with $U$ and are limited in their ability to function as early-warning signals (see figure 3e). This concern has led to considerable doubt in the reliability of CSD-based warning signals for, e.g. climate tipping points [51–53]. However, as we have laid out above, the characteristics of the noise term $\mathrm{d}U$ expressed through its PSD, its autocorrelation and its non-dispersiveness render it a poor choice for the modelling of natural systems. Associated concerns related to the general applicability of early-warning signals for this specific reason are, therefore, exaggerated. In bifurcating systems with stationary red noise in the sense of $U_t\mathrm{d}t$, an increase of variance and AC(1) to infinity and one, respectively, is indeed expected. If noise parameters are non-stationary, however, this may not hold and more specialized statistical methods may be needed to detect CSD [43,54].

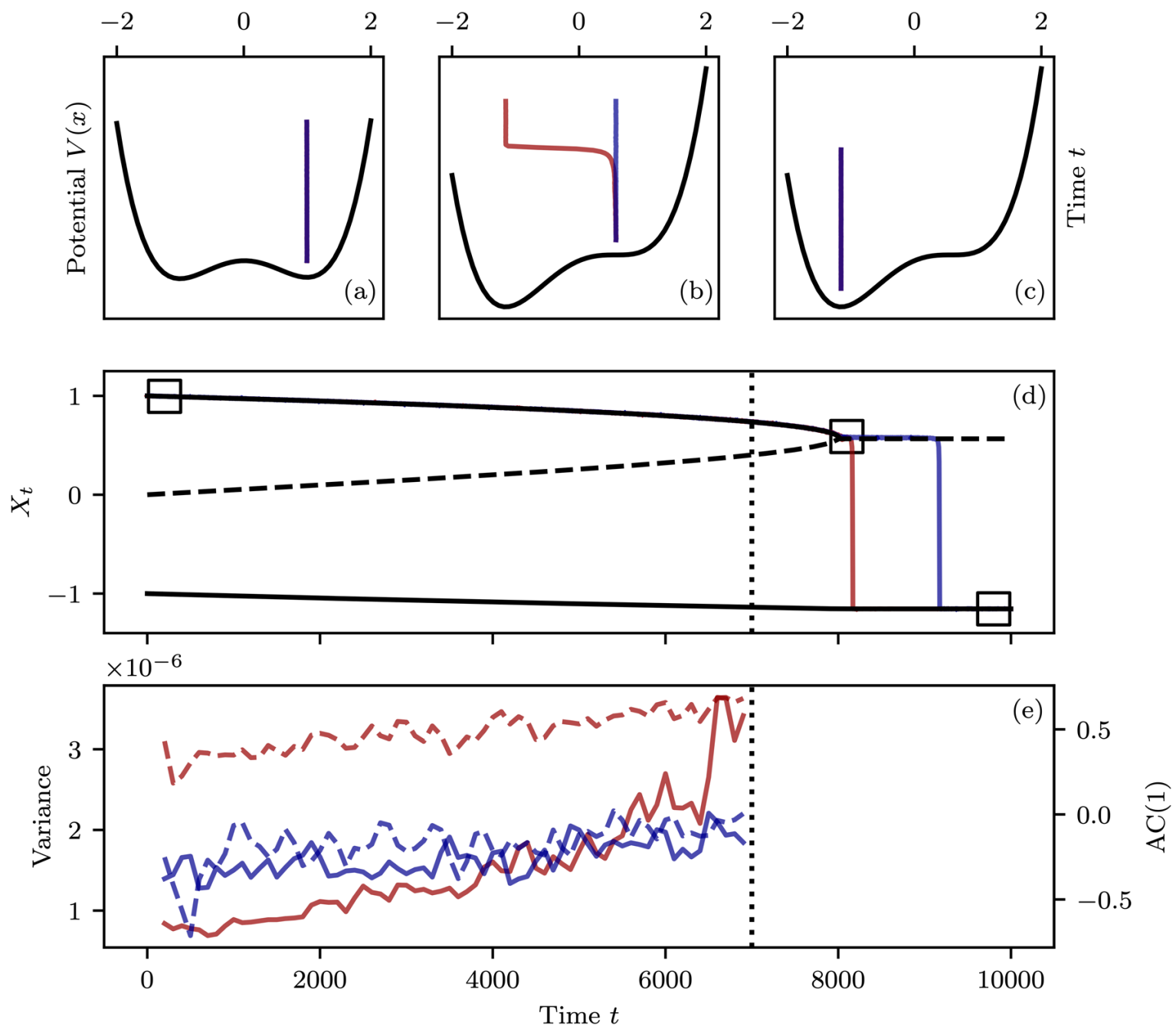

**Figure 3.** Time-series analysis in a conceptual model of a bifurcating system. (a–c) Show a changing potential landscape, first exhibiting two stable equilibria (a) and then one non-hyperbolic and one remaining stable equilibrium (b–c). The time series shown in (a–c) with time on the y-axis are short slices of the time series in panel (d). There, two sample time series of the bifurcating system under different noise forcings are shown. Additionally, the bifurcation diagram of the system is overlaid. The loss of one of the stable states at model time $t = 8000$ incurs an abrupt transition in the system when the system is driven by (natural) dispersive noise. This can be observed for the *red* curves, which correspond to a system driven by red noise $U_t\mathrm{d}t$. The non-dispersive noise model $\mathrm{d}U$ driving the model shown by the *blue* curves incurs a quasi-stable state at the non-hyperbolic equilibrium (the formerly stable equilibrium). It has been shown that in systems driven by $\mathrm{d}U$ noise, the ability to anticipate an approaching tipping point through so-called early-warning signals is diminished. This is because the asymptotic limit of, e.g. system variance is not infinity but a finite quantity [45]. The consequence of this fact can be observed in (e), where the increases in variance (solid line) and lag−1 autocorrelation (*dashed line*) are much less pronounced than for the system driven by red noise $U\mathrm{d}t$. The respective Kendall's $\tau$ statistics quantifying the significance of the increases in variance and AC(1) are $\tau = 0.47$ and 0.43 as opposed to $\tau = 0.87$ and 0.72 in the $U_t\mathrm{d}t$ driven system.

We have seen how the choice between continuous-time red noise $U_t\mathrm{d}t$ and the term $\mathrm{d}U$ can have an impact on the physical aptness of the model to represent natural phenomena. We now turn to another example of where the choice of $\mathrm{d}U$ may cause confusion. This is the case for the two studies [24] and [55] by Bercu *et al.* Both claim to investigate morally the same stochastic system, but in continuous and discrete time, respectively. Statistical time-series estimators for model parameters are derived. In the first study, this is done for the model

$$\mathrm{d}X_t = -\lambda X_t\mathrm{d}t + \mathrm{d}U_t$$
$$\mathrm{d}U_t = -\theta U_t\mathrm{d}t + \mathrm{d}W_t,$$

while the second study investigates

$$X_{k+1} = \lambda X_k + U_k$$
$$U_{k+1} = \theta U_k + z_k,$$

where the $(z_k)_{k\in\mathbb{N}}$ are again i.d.d. unit Gaussian. A serious conceptual issue arises under this supposed model equivalence. The time-series estimators for, e.g. $\lambda$ will each perform well on time series of the given continuous- and discrete-time model, respectively. But the assumption that they would both be applicable to time series from a physical system with a single underlying characteristic is wrong. In the first case, the estimators apply to a physical system driven by negatively correlated noise, while in the second case, the noise is positively correlated in time. The equivocation of $dU$ with discrete-time red noise thus again has tangible consequences for the mathematical analysis of physical systems, a central message of this manuscript.

# 4. Discussion

The decision to restrict ourselves to models of the Itô-differential type is mainly motivated by the prevalence of this framework in the theoretical and applied literature. There are, however, alternatives to this setting which may also allow for sensible red noise models [56–58]. Here, we discuss one notable framework from the literature.

In many applied fields of stochastic modelling, the possibility of using self-similar processes like fractional Brownian motion $B_t^H$ with Hurst parameter $H \neq 1/2$ in lieu of semimartingales for modelling noise has gained great popularity in applications, cf. e.g. [59,60]. The non-stationary autocovariance ($\tau \geq 0$)

$$R_{B^H}(t,\tau) = \frac{1}{2}(t^{2H} + (t+\tau)^{2H} - \tau^{2H})$$

implies stationary increments $B_{t+\Delta t}^H - B_t^H \sim \mathcal{N}(0, \Delta t^{2H})$, which are positively correlated if and only if $H > 1/2$. This so-called long-memory case is interesting to us because it allows for the modelling of persistence in the noise increments. Since the paths of $B^H$ have vanishing $p$-variation for any $p > H^{-1}$, the class of admissible integrand processes $\gamma$ in the Young framework of integration is quite large: if the paths of $\gamma$ possess finite two-variation, then the pathwise Riemann–Stieltjes integral is well-defined. It is, therefore, sensible to consider stochastic dynamical systems in the form of integral equations where some terms are integrated with respect to fractional Brownian motion. Extensive and rigorous theory on this topic may be found in [56]. A consistent derivation of the PSD associated with the differential $dB_t^H$, often referred to as fractional Gaussian noise, is given in [61] (see also figure 1d):

$$S_{dB_t^H}(\omega) = C_H \omega^{1-2H},$$

where $C_H$ is a constant depending only on $H$. Hence, the differential exhibits a vanishing PSD in the limit of high frequencies and may generally be considered an alternative continuous-time modelling approach for the concept of red noise. However, the persistence of fractional Gaussian noise in terms of its autocovariance only decays with $\mathcal{O}(\tau^{2H-2})$, in contrast to the usual exponential decay seen in the Ornstein–Uhlenbeck process (see figure 2d). Fractional Gaussian noise is instead often referred to as a long-memory noise term [62–64]. For red noise models consisting of a combination of an Itô term and a fractional Gaussian noise term, our main result constraining the white noise term continues to hold (see appendix A.2).

Stochastic models of natural systems arise in many applied fields of natural science. The stochastic component of such a model can represent fast, unresolved dynamics of a multi-scale system or express uncertainties in the observed dynamics [1–3]. When constructing stochastic models for natural systems, all available modelling frameworks should be considered based on their conceptual and quantitative ability to represent the dynamics in question (cf. figure 3). Because of its prevalence and precedent in theoretical and applied literature, the Itô-framework of stochastic integration is often the avenue of choice. We have investigated here the implications of certain conceptual notions regarding red noise. In particular, this was the property of a vanishing PSD in infinitely high frequencies, along with practical regularity constraints. The impact of our findings on any concrete modelling choice should be considered carefully. While in some cases, the asymptotic characteristic of the PSD may be explicitly desired, it may not be of much importance in other cases. For such modelling problems, the statements of our Theorem hold no direct insights. We were nevertheless able to discern specifically the differential $dU$ of the Ornstein–Uhlenbeck process as a generally unsuitable model for red noise based on our findings. From a conceptual standpoint, the results and discussion given here should clarify the use of red noise models in continuous-time stochastic modelling.

**Ethics.** This work did not require ethical approval from a human subject or animal welfare committee.

**Data accessibility.** This article has no additional data.

**Declaration of AI use.** We have not used AI-assisted technologies in creating this article.

**Authors' contributions.** A.M.: conceptualization, formal analysis, investigation, methodology, visualization, writing—original draft, writing—review and editing; D.K.: conceptualization, formal analysis, methodology, supervision, writing—original draft, writing—review and editing; N.B.: conceptualization, funding acquisition, project administration, supervision, writing—original draft, writing—review and editing.

All authors gave final approval for publication and agreed to be held accountable for the work performed therein.

**Conflict of interest declaration.** We declare we have no competing interests.

**Funding.** This work has received funding from the Volkswagen Stiftung and under the Marie Sklodowska-Curie grant agreement No. 956170, as well as from the Federal Ministry of Education and Research under grant No. 01LS2001A. This is ClimTip contribution #74; the ClimTip project has received funding from the European Union's Horizon Europe research and innovation programme under grant agreement No. 101137601: Funded by the European Union. Views and opinions expressed are, however, those of the author(s) only and do not necessarily reflect those of the European Union or the European Climate, Infrastructure and Environment Executive Agency (CINEMA). Neither the European Union nor the granting authority can be held responsible for them.

# Appendix A.

## A.1. Proof of theorem 2.1

1. *Finite time horizon*. First, we show that the constant term in the finite-time PSD stemming from the $\beta_t \mathrm{d}W_t$ term must be matched by the PSD of $\alpha_t \mathrm{d}t$ to result in an overall PSD which is vanishing at $\omega \to \infty$.

$$
\begin{aligned}
S^{(T)}_{\mathrm{d}Y}(\omega) &= \mathbb{E}\left[\frac{1}{T}\left|\int_0^T \exp(-i\omega t)\mathrm{d}Y_t\right|^2\right] \\
&= \mathbb{E}\left[\frac{1}{T}\left|\int_0^T \exp(-i\omega t)\alpha_t \mathrm{d}t + \int_0^T \exp(-i\omega t)\beta_t \mathrm{d}W_t\right|^2\right] \\
&= \frac{1}{T}\left\|\int_0^T \exp(-i\omega t)\alpha_t \mathrm{d}t + \int_0^T \exp(-i\omega t)\beta_t \mathrm{d}W_t\right\|^2_{L^2(\Omega)} \\
&\geq \left[\frac{1}{\sqrt{T}}\left\|\int_0^T \exp(-i\omega t)\alpha_t \mathrm{d}t\right\|_{L^2(\Omega)} - \frac{1}{\sqrt{T}}\left\|\int_0^T \exp(-i\omega t)\beta_t \mathrm{d}W_t\right\|_{L^2(\Omega)}\right]^2 \\
&= \left(\sqrt{S^{(T)}_{\alpha_t \mathrm{d}t}(\omega)} - \sqrt{\frac{1}{T}\int_0^T \mathbb{E}\left[\beta_t^2\right]\mathrm{d}t}\right)^2,
\end{aligned}
$$

where we used the inverse triangle inequality. Since we assumed $S^{(T)}_{\mathrm{d}Y}(\omega) \xrightarrow[\omega\to\infty]{} 0$, we conclude

$$
S^{(T)}_{\alpha_t \mathrm{d}t}(\omega) \xrightarrow[\omega\to\infty]{} \frac{1}{T}\int_0^T \mathbb{E}\left[\beta_t^2\right]\mathrm{d}t.
$$

It remains to be proven that this limit is 0. Since we have uniform integrability of the random variables $\int_0^T \exp(-i\omega t)\alpha_t \mathrm{d}t$ for all choices of $\omega \in \mathbb{R}$ through

$$
\left\|\int_0^T \exp(-i\omega t)\alpha_t \mathrm{d}t\right\|_{L^2(\Omega)} \leq \int_0^T \|\alpha_t\|_{L^2(\Omega)}\,\mathrm{d}t < \infty,
$$

it suffices to prove that

$$
\int_0^T \exp(-i\omega t)\alpha_t \mathrm{d}t \xrightarrow[\omega\to\infty]{a.s.} 0. \tag{A 1}
$$

But each path $a : [0, T] \to \mathbb{R}$ of $\alpha$ is uniformly continuous on the compact domain, i.e. for any $\epsilon > 0$, we find a $\delta_\epsilon > 0$ such that $|t_1 - t_0| \leq \delta_\epsilon \Rightarrow |a(t_1) - a(t_0)| \leq \epsilon$. Define $c := \max_{t\in[0,T]} a(t)$ and consider

$$\begin{aligned}\left|\int_{t_0}^{t_0+\delta_\epsilon} \cos(\omega t)a(t)\mathrm{d}t\right| &= \left|\int_{t_0}^{t_0+\delta_\epsilon} \cos(\omega t)a(t_0)\mathrm{d}t + \int_{t_0}^{t_0+\delta_\epsilon} \cos(\omega t)(a(t) - a(t_0))\mathrm{d}t\right| \\ &\leq c\left|\int_{t_0}^{t_0+\delta_\epsilon} \cos(\omega t)\mathrm{d}t\right| + \epsilon\int_{t_0}^{t_0+\delta_\epsilon} |\cos(\omega t)|\,\mathrm{d}t \\ &\leq c\frac{2}{\omega} + \epsilon\delta_\epsilon.\end{aligned}$$

Piecing together $\lceil T/\delta_\epsilon\rceil$ intervals of length of at most $\delta_\epsilon$, we obtain for frequencies $\omega \geq \omega_\epsilon := \frac{1}{\delta_\epsilon \epsilon}$

$$\left|\int_0^T \cos(\omega t)a(t)\mathrm{d}t\right| \leq \left(c\frac{2}{\omega} + \epsilon\delta_\epsilon\right)\left\lceil\frac{T}{\delta_\epsilon}\right\rceil \leq (2c\epsilon + \epsilon)\,2T =: \kappa_\epsilon.$$

This implies that for sufficiently large frequencies, the left-hand side can be bounded by any value $\kappa_\epsilon > 0$. The limit of the left-hand side as $\omega \to \infty$ must therefore be 0.

The real part and, similarly, the imaginary part of (A 1) converge to 0 for every path of $\alpha$. Hence, the limit in question is also 0:

$$S^{(T)}_{\alpha_t \mathrm{d}t}(\omega) \xrightarrow[\omega\to\infty]{} 0 = \frac{1}{T}\int_0^T \mathbb{E}\left[\beta_t^2\right]\mathrm{d}t.$$

We may conclude that for almost all $t \in [0, T]$ we have $\beta_t = 0$ $\mathbb{P}$-almost surely.

2. *Infinite time horizon*. By identical arguments as in the first part of the proof, we arrive at

$$S_{\alpha_t \mathrm{d}t}(\omega) := \lim_{T\to\infty} \mathbb{E}\left[\frac{1}{T}\left|\int_0^T \exp(-i\omega t)\alpha_t \mathrm{d}t\right|^2\right] \xrightarrow[\omega\to\infty]{} \lim_{T\to\infty}\frac{1}{T}\int_0^T \mathbb{E}\left[\beta_t^2\right]\mathrm{d}t = \mathbb{E}\left[\beta_0^2\right],$$

on the condition that the limit of $T \to \infty$ in the above definition of $S_{\alpha_t \mathrm{d}t}$ exists. Because $\alpha$ is square-integrable, its autocovariance $R_\alpha(\tau)$ is bounded by $\mathrm{Var}(\alpha_0)$ for all $\tau \in \mathbb{R}$. This, together with the absolute integrability of $R_\alpha$ implies that $R_\alpha \in L^2(\mathbb{R})$. By the Wiener-Khinchin theorem and Plancherel's theorem, we deduce

$$S_{\alpha_t \mathrm{d}t} = \mathcal{F}[R_\alpha] \in L^2(\mathbb{R})$$

and thereby confirm that the limit in the definition of $S_{\alpha_t \mathrm{d}t}$ exists. We have established that $S_{\alpha_t \mathrm{d}t}$ converges to a constant as $\omega \to \infty$. In order for $S_{\alpha_t \mathrm{d}t}$ to be square-integrable, the constant in question must be 0. Hence, we must have

$$\mathbb{E}\left[\beta_t^2\right] = \mathbb{E}\left[\beta_0^2\right] = 0$$

which implies $\beta_t = 0$ $\mathbb{P}$-a.s. for all $t \geq 0$. ■

## A.2. Extension of the theorem for fractional Gaussian noise terms

Some models of physical systems use a fractional Gaussian noise term $\mathrm{d}B_t^H$ with $H > 1/2$ to represent dynamics exhibiting (long) memory in time. As mentioned in §4, the PSD of such a term vanishes in the limit of large frequencies

$$S_{\mathrm{d}B_t^H}(\omega) = C_H \omega^{1-2H},$$

making the term a viable candidate to model red noise. In this section, we show that even in the presence of a fractional Gaussian noise term $\gamma(t)\mathrm{d}B_t^H$, any red noise model

$$\mathrm{d}Y_t = \alpha_t \mathrm{d}t + \beta_t \mathrm{d}W_t + \gamma(t)\mathrm{d}B_t^H \tag{A 2}$$

is similarly restricted in the choice of $\beta$.

**Lemma A.1.** *Let $\alpha$ and $\beta$ be adapted processes satisfying the integrability condition (2.1). Suppose that $\alpha$ has continuous paths, $\beta$ is predictable, and the function $\gamma:\mathbb{R}_+\to\mathbb{R}$ is of bounded variation and Hölder-continuous with exponent $1-H<\alpha_\gamma\le 1$. Define the stochastic process*

$$Y_t=\int_0^t \alpha_s \mathrm{d}s+\int_0^t \beta_s \mathrm{d}W_s+\int_0^t \gamma(s)\mathrm{d}B_s^H,\quad t\ge 0,$$

*where $B_t^H$ is a fractional Brownian motion with Hurst index $H>1/2$, adapted to $(\mathcal{F}_t)_{t\in\mathbb{R}_+}$.*

*1. Finite time horizon. Assume that the finite-time PSD of $\mathrm{d}Y$ vanishes in the limit of infinitely high frequencies, i.e.*

$$S_{\mathrm{d}Y}^{(T)}(\omega):=\mathbb{E}\left[\frac{1}{T}\left|\int_0^T \exp(-i\omega t)\mathrm{d}Y_t\right|^2\right]\xrightarrow[\omega\to\infty]{}0.$$

*Then $\beta_t=0$ $\mathbb{P}$-a.s. for almost all $t\in[0,T]$.*

*2. Infinite time horizon. Assume that $\alpha$ and $\beta$ are stationary processes and let $\alpha$ be centred around 0 with an absolutely integrable autocovariance structure $R_\alpha(\tau)$. Let $\gamma\equiv\gamma_0\in\mathbb{R}$ be a constant. Assume that the PSD of $\mathrm{d}Y$ on an infinite time horizon exists and that it vanishes in the limit of infinitely high frequencies, i.e.*

$$S_{\mathrm{d}Y}(\omega):=\lim_{T\to\infty}\mathbb{E}\left[\frac{1}{T}\left|\int_0^T \exp(-i\omega t)\mathrm{d}Y_t\right|^2\right]\xrightarrow[\omega\to\infty]{}0.$$

*Then $\beta_t=0$ $\mathbb{P}$-a.s. for all $t\ge 0$.*

*Proof.* The line of argumentation is identical to that put forth in the proof of the main result. We need only prove

$$S^{(T)}_{\gamma(t)\mathrm{d}B_t^H}(\omega)\xrightarrow[\omega\to\infty]{}0$$
$$\text{and}\quad S_{\gamma(t)\mathrm{d}B_t^H}(\omega)\xrightarrow[\omega\to\infty]{}0$$

for the finite and infinite time case, respectively.

1. *Finite time horizon.* We can apply partial integration to obtain

$$\begin{aligned}\left|\int_0^T \exp(-i\omega t)\gamma(t)\mathrm{d}B_t^H\right|^2 &= \left|\gamma(T)\int_0^T \exp(-i\omega t)\mathrm{d}B_t^H-\int_0^T\int_0^s \exp(-i\omega t)\mathrm{d}B_t^H\mathrm{d}\gamma(s)\right|^2\\ &= \left|\gamma(T)\int_0^T \exp(-i\omega t)\mathrm{d}B_t^H\right|^2+\left|\int_0^T\int_0^s \exp(-i\omega t)\mathrm{d}B_t^H\mathrm{d}\gamma(s)\right|^2\\ &\quad+\left|\gamma(T)\int_0^T \exp(-i\omega t)\mathrm{d}B_t^H\int_0^T\int_0^s \exp(-i\omega t)\mathrm{d}B_t^H\mathrm{d}\gamma(s)\right|.\end{aligned}\tag{A 3}$$

All integrals are well-defined within the Young framework since we have imposed an appropriate constraint on the Hölder-exponent of $\gamma$. The first term taken in expectation will tend to 0 in the limit of large frequencies. The second term can be bounded by introducing the monotonously increasing function $\gamma^+:=\gamma_\uparrow-\gamma_\downarrow$, where the latter is defined as $\gamma$ restricted to the domain where it is increasing and decreasing, respectively. This is possible since $\gamma$ is assumed to be of bounded variation.

$$\left|\int_0^T\int_0^s \exp(-i\omega t)\mathrm{d}B_t^H\mathrm{d}\gamma(s)\right|^2\le\int_0^T\left|\int_0^s \exp(-i\omega t)\mathrm{d}B_t^H\right|^2\mathrm{d}\gamma^+(s)$$

Taking the expectation of the inner integral, we see that this term, too, vanishes in the limit of infinitely large frequencies. The last term of (B 2) may be treated similarly.

2. *Infinite time horizon.* Since it is assumed that $\gamma\equiv\gamma_0\in\mathbb{R}$ is a constant in this case, we may conclude the proof with the observation that

$$S_{\gamma_0\mathrm{d}B_t^H}(\omega)=\gamma_0^2\,S_{\mathrm{d}B_t^H}(\omega)\xrightarrow[\omega\to\infty]{}0.$$

■

## A.3. Derivation of the power spectral densities

For unit white noise d$W$, the calculation is a straightforward application of the Itô-isometry:

$$S_{\mathrm{d}W}(\omega) = \lim_{T\to\infty} \mathbb{E}\left[\frac{1}{T}\left|\int_0^T \exp(-i\omega t)\mathrm{d}W_t\right|^2\right] = \lim_{T\to\infty} \mathbb{E}\left[\frac{1}{T}\int_0^T 1\mathrm{d}t\right] = 1.$$

For the other claims, first define the Ornstein–Uhlenbeck process through the SDE

$$\mathrm{d}U_t = -\theta U_t\mathrm{d}t + \mathrm{d}W_t, \quad U_0 \sim \mathcal{N}\left(0, \frac{1}{2\theta}\right),$$

so that $U$ is stationary with autocovariance $\mathrm{Cov}\,(U_t, U_{t+\tau}) = \frac{1}{2\theta}\exp(-\theta\,|\tau|)$ and compute via the Fubini–Tonelli Theorem

$$\begin{aligned}
S_{U_t\mathrm{d}t}(\omega) &= \lim_{T\to\infty} \mathbb{E}\left[\frac{1}{T}\left|\int_0^T \exp(-i\omega t)U_t\mathrm{d}t\right|^2\right] \\
&= \lim_{T\to\infty}\frac{1}{T}\int_0^T\int_0^T \exp(-i\omega(t-s))\mathrm{Cov}\,(U_t, U_s)\,\mathrm{d}t\mathrm{d}s \\
&= \lim_{T\to\infty}\frac{1}{T}\int_0^T\int_0^T \exp(-i\omega(t-s))\frac{1}{2\theta}\exp(-\theta\,|t-s|)\mathrm{d}t\mathrm{d}s \\
&= \lim_{T\to\infty}\frac{1}{T}f_1(T,\omega,\theta) = \frac{1}{\theta^2+\omega^2}.
\end{aligned}$$

The function $f_1$ can be computed to be

$$f_1(T,\omega,\theta) = \frac{T}{\theta^2+\omega^2} + \frac{(\omega^2-\theta^2)(1-\exp(-\theta T)\cos(\omega T)) - 4\theta\omega\exp(-\theta T))\sin(T\omega)}{\theta(\theta^2+\omega^2)}$$

and we will make use of it again in the derivation of $S_{\mathrm{d}U}$:

$$\begin{aligned}
S_{\mathrm{d}U}(\omega) &= \lim_{T\to\infty} \mathbb{E}\left[\frac{1}{T}\left|\int_0^T \exp(-i\omega t)\mathrm{d}U_t\right|^2\right] \\
&= \lim_{T\to\infty} \mathbb{E}\left[\frac{1}{T}\left|\int_0^T \exp(-i\omega t)(-\theta)U_t\mathrm{d}t + \int_0^T \exp(-i\omega t)\mathrm{d}W_t\right|^2\right] \\
&= \lim_{T\to\infty}\frac{1}{T}\left[\theta^2 f_1(T,\omega,\theta) - \theta\mathbb{E}\left[\int_0^T \exp(-i\omega t)\mathrm{d}W_t\int_0^T \exp(i\omega t)U_t\mathrm{d}t\right]\right. \\
&\qquad \left. - \theta\mathbb{E}\left[\int_0^T \exp(-i\omega t)U_t\mathrm{d}t\int_0^T \exp(i\omega t)\mathrm{d}W_t\right] + T\right] \\
&= \lim_{T\to\infty}\frac{1}{T}\left[\theta^2 f_1(T,\omega,\theta) - \theta f_2(T,\omega,\theta) - \theta\overline{f_2(T,\omega,\theta)} + T\right]
\end{aligned}$$

where we introduced

$$\begin{aligned}
f_2(T,\omega,\theta) &:= \mathbb{E}\left[\int_0^T \exp(-i\omega t)\mathrm{d}W_t\int_0^T \exp(i\omega t)U_t\mathrm{d}t\right] \\
&= \mathbb{E}\left[\int_0^T\int_0^T \exp(-i\omega(s-t))U_t\mathrm{d}W_s\mathrm{d}t\right] \\
&= \mathbb{E}\left[\int_0^T\int_0^T\int_0^t \exp(-i\omega(s-t))\exp(-\theta(t-r))\mathrm{d}W_r\mathrm{d}W_s\mathrm{d}t\right].
\end{aligned}$$

In the last step, we inserted the closed-form solution to the Ornstein–Uhlenbeck SDE, where $U_0$ is independent of $W_t$ for all $t \geq 0$:

$$U_t = U_0\exp(-\theta t) + \int_0^t \exp(-\theta(t-r))\mathrm{d}W_r.$$

The integration bounds in the integral with respect to $\mathrm{d}W_s$ may as well be $(0,t)$ since the integral in $(t,T)$ is independent of the innermost integral, and their expectations are 0 each. Using the Itô isometry again, we get

$$f_2(T,\omega,\theta) = \int_0^T \int_0^t \exp(-i\omega(s-t))\exp(-\theta(t-s))\mathrm{d}s\mathrm{d}t$$
$$= \frac{T(\theta-i\omega)+e^{-T(\theta-i\omega)}-1}{(\theta-i\omega)^2}.$$

We see that

$$\lim_{T\to\infty}\frac{1}{T}f_1(T,\omega,\theta) = \frac{1}{\theta^2+\omega^2} = \lim_{T\to\infty}\frac{1}{T}\frac{1}{2\theta}\left(f_2(T,\omega,\theta)+\overline{f_2(T,\omega,\theta)}\right)$$

and so the spectral density is simply

$$S_{\mathrm{d}Y}(\omega) = \lim_{T\to\infty}\frac{1}{T}\left[\theta^2 f_1(T,\omega,\theta)-\theta f_2(T,\omega,\theta)-\theta\overline{f_2(T,\omega,\theta)}+T\right]$$
$$= \frac{-\theta^2}{\theta^2+\omega^2}+1 = \frac{\omega^2}{\theta^2+\omega^2}.$$